\documentclass{amsproc}
\usepackage{amsfonts}
\usepackage{amsthm}

\theoremstyle{plain}

\newtheorem{corollary}{Corollary}

\newtheorem{example}{Example}

\newtheorem{remark}{Remark}

\newtheorem*{theorem}{Theorem}
\input{tcilatex}

\begin{document}
\title{Nonlinear Convergence Sets of Divergent Power Series}
\author{Buma L.~Fridman, Daowei Ma and Tejinder S. Neelon}
\dedicatory{To Professor J\'{o}zef Siciak on his 80th birthday}

\begin{abstract}
A nonlinear generalization of convergence sets of formal power series, in
the sense of Abhyankar-Moh\cite{AM}, is introduced. Given a family $%
y=\varphi _{s}(t,x)=sb_{1}(x)t+b_{2}(x)t^{2}+...$ of analytic curves in $%
\mathbb{C\times C}^{n}$ passing through the origin, $\limfunc{Conv}%
\nolimits_{\varphi }(f)$ of a formal power series $f(y,t,x)\in \mathbb{C}%
[[y,t,x]]$ is defined to be the set of all $s\in \mathbb{C}$ for which the
power series $f(\varphi _{s}(t,x),t,x)$ converges as a series in $(t,x).$ We
prove that for a subset $E\subset \mathbb{C}$ there exists a \emph{divergent}
formal power series $f(y,t,x)\in \mathbb{C}[[y,t,x]]$ such that $E=\limfunc{%
Conv}_{\varphi }(f)$ if and only if $E$ is a $F_{\sigma }$ set of zero
capacity. This generalizes the results of P. Lelong and A. Sathaye for the
linear case $\varphi _{s}(t,x)=st.$
\end{abstract}

\keywords{convergence sets, formal power series, analytic functions,
capacity, polar sets.}
\subjclass[2000]{Primary: 26E05, 32A05, 30C85, 40A05}
\address{ fridman@math.wichita.edu, Department of Mathematics, Wichita State
University, Wichita, KS 67260-0033, USA}
\address{ dma@math.wichita.edu, Department of Mathematics, Wichita State
University, Wichita, KS 67260-0033, USA}
\address{ neelon@csusm.edu, Department of Mathematics, California State
University San Marcos, San Marcos CA 92096-0001, USA}
\maketitle

\bigskip We say that formal power series $f(z)=\sum_{|\alpha |=0}^{\infty
}a_{\alpha }z^{\alpha }$ is convergent if there exists a constant $C$ such
that $|a_{\alpha }|\leq C^{|\alpha |},$ for all $\alpha \in \mathbb{Z}%
_{+}^{n}$. (Here we have used multiindex notation: $\mathbb{Z}_{+}^{n}$
denote the set of all $n$-tuples $\alpha :=(\alpha _{1},\alpha
_{2},...,\alpha _{n})$ of integers $\alpha _{i}\geq 0,$ if $%
z=(z_{1},z_{2},...,z_{n})$ and $\alpha \in \mathbb{Z}_{+}^{n},$ then $%
z^{\alpha }=$ $z_{1}^{\alpha _{1}}z_{2}^{\alpha _{2}}...z_{n}^{\alpha _{n}},$
and $|\alpha |:=\alpha _{1}+\alpha _{2}+...+\alpha _{n}$ denotes the length
of $\alpha \in \mathbb{Z}_{+}^{n}.$) A series $f$ is called divergent if it
is not convergent. A divergent power series may still converge when
restricted to a certain set of lines or planes through the origin. For
example, Abhyankar and Moh \cite{AM} considered the convergence set $%
\limfunc{Conv}(f)$ of a series $f$ defined to be the set of all $s\in 
\mathbb{C}$ for which $f(sz_{2},z_{2},...,z_{n})$ converges as a series in $%
(z_{2},z_{3},...,z_{n}).$ The convergence set of divergent series can be
empty or an arbitrary countable set (see examples below). The Abhyankar-Moh
paper proved that the one dimensional Hausdorff measure of the convergence
set of a divergent series is zero. In the case when $n=2,$ Pierre Lelong had
earlier proved that if $\limfunc{Conv}(f)$ is not contained in a $F_{\sigma
} $ set of zero capacity then the series $f$ is necessarily convergent, and
conversely, given any set \ $E$\ contained in a $F_{\sigma }$\ set of zero
capacity a divergent power series $f$ can be constructed so that $E\subset 
\limfunc{Conv}(f).$ \ This result has been rediscovered, independently, by
several authors (see e.g. \cite{LM}, \cite{Ne}, \cite{Sa}, and see also \cite%
{FM}, \cite{FM1}, \cite{Ne1},\cite{Ri}, for other related results). The
optimal result was obtained by Sathaye \cite{Sa} who strengthened the
results of Abhyankar-Moh and Lelong by proving that a necessary and
sufficient condition for a set $S\subset C$ to be equal to the convergence
set of a divergent power series $f(z)$\ \ is that $S$\ is $F_{\sigma }$-set
of transfinite diameter zero i.e. $S=\cup _{j=1}^{\infty }E_{j}$\ where each 
$E_{j}$\ is closed set of transfinite diameter zero. These results can be
viewed as an optimal and formal analogs of Hartogs' Theorem on separate
analyticity in several complex variables.

In this article, we consider `nonlinear' convergence sets of formal power
series $f(y,t,x)\in \mathbb{C}[[y,t,x]],$ $x=(x_{1},x_{2},...,x_{n})$ by
restricting $f(y,t,x)$ along a one-parameter family of `tangential'
perturbations of a fixed analytic curve $y=$ $\varphi (t,x)$ through the
origin.\textrm{\ }Throughout this paper, $\varphi (t,x):=\sum_{j=1}^{\infty
}b_{j}(x)t^{j}$ will denote a fixed convergent power series where $%
b_{j}(x):=\sum_{i=0}^{\infty }b_{ji}x^{i},j=1,2,3,...,$ are convergent power
series in $x$ with complex coefficients. We assume that $b_{10}=1$.

For $s\in \mathbb{C},$ we put $\varphi _{s}(t,x)=\varphi
(s,t,x)=sb_{1}(x)t+\sum_{j=2}^{\infty }b_{j}(x)t^{j}$. Define the $\varphi $%
-convergence set of a series $f(y,t,x)\in \mathbb{C}[[y,t,x]],$ as follows. 
\begin{equation*}
\limfunc{Conv}\nolimits_{\varphi }(f):=\{s\in \mathbb{C}:f(\varphi
(s,t,x),t,x)\;\text{ converges as a series in }(t,x)\}.
\end{equation*}

\ If $K$ is an infinite compact set\ in $\mathbb{C},$ then for each positive
integer $n$ there is a unique monic polynomial $P_{n,K}(z)$ (called the
Chebychev polynomial)\ such that \ 
\begin{equation*}
\rho _{n}(K):=\max_{z\in K}\left\vert P_{n,K}(z)\right\vert
=\min_{P_{n}}(\max_{z\in K}|P_{n}(z)|)\text{ }
\end{equation*}%
where minimum is taken over the set of all monic polynomials $P_{n}(z)=$\ $%
z^{n}+a_{1}z^{n-1}+\cdots +a_{n}$. The limit $\rho (K):=\lim_{n}\left( \rho
_{n}(K)\right) ^{1/n}$\ exists and is called the Chebychev constant of $K.$
For a compact subset $K\ $of $\mathbb{C},$ $\rho (K)\ $coincides with the
logarithmic capacity $\mathfrak{c}(K)$ and the transfinite diameter $%
d_{\infty }(K)$\ of $K$\ (see \cite{Ah}, Chapter 2). If \ $E=\cup
_{n=1}^{\infty }K_{n}$, where $K_{n}$ are compact sets of zero capacity,
then $\mathfrak{c}(E)=0.$ A subset $E$ of $\mathbb{C}$ is zero capacity if
and only if it is polar i.e. $E\subset \{z:u(z)=-\infty \}$ for some
nonconstant subharmonic function $u:\mathbb{C}\rightarrow \lbrack -\infty
,\infty )$. A $F_{\sigma }$ set $E$ in $\mathbb{C}$ is said to have zero
capacity if $\mathfrak{c}(E\cap \{|z|\leq r\})=0$ for every $r\geq 0.$

\begin{theorem}
Let $\varphi (t,x)$ be as above, and let $E$ be a subset of $\mathbb{C}$.
There exists a \emph{divergent} formal power series $f(y,t,x)\in \mathbb{C}%
[[y,t,x]]$ such that $E=\limfunc{Conv}_{\varphi }(f)$ if and only if $E$ is
an $F_{\sigma }$ set of zero capacity.
\end{theorem}

\begin{proof}
Suppose that $E$ is an $F_{\sigma }$ set with $\mathfrak{c}(E)>0$. By
replacing $E$ with a compact subset \ $K\subset E$ of positive capacity we
can assume that $E$ is compact. Let $f(y,t,x):=%
\sum_{i,j,k}a_{ijk}y^{i}t^{j}x^{k}\in \mathbb{C[}[y,t,x]]$ be such that $%
g(s;t,x):=f(\varphi (s,t,x),t,x)$ converges for each $s\in E$. We need to
show that $f$ is convergent. Rewrite $f$ as $f(y,t,x)=%
\sum_{i,j}a_{ij}(x)y^{i}t^{j},$ where $a_{ij}(x):=\sum_{k=0}^{\infty
}a_{ijk}x^{k}\in \mathbb{C}[[x]],$ and 
\begin{equation*}
g(s;t,x):=\sum_{q\geq 0,k\geq 0}\lambda _{qk}(s)t^{q}x^{k}:=\sum_{p\geq
0,q\geq 0}d_{pq}(x)s^{p}t^{q}.
\end{equation*}

It is clear that $\lambda _{qk}(s)$ is a polynomial of degree at most $q$,
and thus 
\begin{equation*}
d_{pq}(x)=0,\;\;\text{for }p>q.
\end{equation*}%
Let $d_{pq}(x):=\sum_{k\geq 0}d_{pqk}x^{k}$, and write $\lambda
_{qk}(s):=\sum_{p=0}^{q}d_{pqk}s^{p}$.

We have 
\begin{equation*}
d_{pq}(x)={\sum\nolimits}^{\prime}a_{ij}(x){\frac{i!}{p!m_{2}!m_{3}!\cdots }}%
b_{1}(x)^{p}b_{2}(x)^{m_{2}}b_{3}(x)^{m_{3}}\cdots ,
\end{equation*}%
where the summation $\sum^{\prime}$ is taken over all nonnegative integers $%
i,j,m_{2},m_{3},\dots $ satisfying 
\begin{equation*}
j+p+2m_{2}+3m_{3}+\cdots =q\text{ and }p+m_{2}+m_{3}+\cdots =i.
\end{equation*}%
Since 
\begin{equation}
\begin{array}{lll}
d_{qq}(x) & = & a_{q,0}(x)b_{1}(x)^{q}, \\ 
d_{q-1,q}(x) & = & a_{q-1,1}(x)b_{1}(x)^{q-1}, \\ 
d_{q-2,q}(x) & = & 
a_{q-2,2}(x)b_{1}(x)^{q-2}+a_{q-1,0}(x)(q-1)b_{1}(x)^{q-2}b_{2}(x), \\ 
d_{q-k,q}(x) & = & a_{q-k,k}(x)b_{1}(x)^{q-k}+\text{terms involving }%
a_{ij}(x)\text{ with }i+j<q\text{,}%
\end{array}
\label{exp}
\end{equation}%
it follows that $a_{ij}(x)$ can be solved uniquely in terms of $d_{pq}(x)$.
In particular, if $d_{pq}(x)=0$ for all $p,q$ then all $a_{ij}(x)=0$ for all 
$i,j$.

For each $s\in E$, there is a constant $C_{s}$ such that $|\lambda
_{qk}(s)|\leq C_{s}^{q+k}$ for all $q+k\geq 1$, since the power series $%
\sum_{q,k}\lambda _{qk}(s)t^{q}x^{k}$ converges. For each positive integer $%
n $, set 
\begin{equation*}
E_{n}=\{s\in E:|\lambda _{qk}(s)|\leq n^{q+k}\;\;\forall q+k\geq 1\}.
\end{equation*}%
The sets $E_{n}$ are closed and $E=\cup _{n=1}^{\infty }E_{n}$. There is a
positive integer $N$ such that $E^{\prime }:=\cup _{n=1}^{N}E_{n}$ has
positive capacity. It follows that $|\lambda _{qk}(s)|\leq N^{q+k}$ for $%
q+k\geq 1$ and for $s\in E^{\prime }$. By the Bernstein-Walsh inequality
(see \cite{FM}, Lemma~1.4), there is a constant $C_{E^{\prime }}\geq 1$ such
that $|d_{pqk}|\leq C_{E^{\prime }}^{q}N^{q+k}\leq (C_{E^{\prime }}N)^{q+k}$.

For some $\tau >0$, $g(s;t,x)$ represents a holomorphic function in $\Delta
_{\tau }\times \Delta _{\tau }\times \Delta _{\tau }$, where $\Delta _{\tau
}=\{z\in \mathbb{C}:|z|<\tau \}$. Shrinking $\tau ,$ if necessary, we may
assume that 
\begin{equation*}
\min \{|b_{1}(x)|:x\in \mathbb{C},|x|\leq \tau \}\geq 1/2
\end{equation*}%
and 
\begin{equation}
\sum_{q,k}|b_{qk}|\tau ^{q+k}<\infty ,\;\;\sum_{p,q,k}|d_{pqk}|\tau
^{p+k}2^{q}(\tau +\sum |b_{ij}|\tau ^{i+j-1})^{q}<\infty .  \label{abs}
\end{equation}
The map $\psi :\mathbb{C}^{n+2}\rightarrow \mathbb{C}^{n+2}$ defined by $%
\psi (s,t,x):=(\varphi (s,t,x),t,x)$ is holomorphic near the origin and is
injective on $Q=\{(s,t,x)\in \Delta _{\tau }\times \Delta _{\tau }\times
\Delta _{\tau }:t\not=0\}$. It follows that there is a holomorphic function $%
G(u,v,w)$ defined on $\psi (Q)$ such that $g=G\circ \psi $ on $Q$.

We now prove that $G$ extends holomorphically to a neighborhood of the
origin. Choose a $\delta ,$ $0<\delta <\tau /2$, sufficiently small so that
the set 
\begin{equation*}
\Gamma :=\{(u,v,w)\in \mathbb{C}^{3}:|u|\leq \delta ^{2},|v|=\delta ,|w|\leq
\delta \}
\end{equation*}%
is contained in $\psi (Q).$ The function $G$ extends holomorphically to a
neighborhood of the origin, if for $|u_{0}|<\delta ^{2}$, $|w_{0}|<\delta $, 
\begin{equation*}
I_{k}(u_{0},w_{0}):={\frac{1}{2\pi \sqrt{-1}}}\int_{|v|=\delta
}v^{k}G(u_{0},v,w_{0})\,dv=0,\forall k=0,1,2,\dots .
\end{equation*}%
For fixed $w_{0}$ and $u_{0}$, write 
\begin{equation*}
\left( {\frac{u_{0}-(\varphi (t,w_{0})-tb_1(w_0))}{tb_{1}(w_{0})}}\right)
^{p}:=\sum_{j=-p}^{\infty }c_{pj}(u_{0},w_{0})t^{j}.
\end{equation*}%
By making use of (\ref{abs}), and substituting the above series expansion
into the integrand in $I_{k}(u_{0},w_{0})$, we obtain 
\begin{eqnarray*}
I_{k}(u_{0},w_{0})&=&{\frac{1}{2\pi \sqrt{-1}}}\int_{|t|=\delta
}t^{k}g\left( {\frac{u_{0}-(\varphi (t,w_{0})-tb_1(w_0))}{tb_{1}(w_{0})}}%
,t,w_{0}\right) dt \\
&=&\sum d_{pq}(w_{0})c_{pj}(u_{0},w_{0}),
\end{eqnarray*}%
where the sum is over all $p,q,j$ with $q+j=-k-1$ and $j\geq -p$. Since $%
q+j=-k-1$ and $j\geq -p$ imply $q=-j-k-1<-j\leq p$, and since $d_{pq}=0$ for 
$q<p$, we see that $I_{k}(u_{0},w_{0})=0$ for $k=0,1,2,\dots $ and for all $%
(u_{0},w_{0})$ with $|u_{0}|<\delta ^{2}$, $|w_{0}|<\delta $. Therefore, $G$
extends holomorphically to a neighborhood of the origin.

Now $g=G\circ \psi $ on $\Delta _{\tau }\times \Delta _{\tau }\times \Delta
_{\tau }$. Hence $g=G\circ \psi $ as a formal power series. Since $g=f\circ
\psi $, we see that $\hat{f}\circ \psi :=\sum \hat{d}_{pq}(x)s^{p}t^{q}=0$,
where $\hat{f}:=\sum_{i,j}\hat{a}_{ij}(x)y^{i}t^{j}=f-G$. It follows that
all $\hat{d}_{pq}(x)$, and hence $\hat{a}_{ij}(x)$ are all $0$. This proves $%
f$ is convergent as $\hat{f}=0$ and $f\equiv G$.

Conversely, suppose $E$ is an $F_{\sigma }$ set with $\mathfrak{c}(E)=0$. We
construct a divergent power series $f(y,t,x)$ such that $\limfunc{Conv}%
_{\varphi }(f)=E$.

By Theorem 6.1 of \cite{Ri} , there exists an increasing sequence $\{q_{j}\}$
of positive integers and a sequence of polynomials $\{P_{j}(s)\}$ with $%
\hbox{\rm deg}(P_{j})\leq $ $q_{j}$, for all $j=1,2,3,...,$ such that the
series $\psi _{s}(t)=\sum_{j}P_{j}(s)t^{q_{j}}$ converges for each $s\in E$,
and diverges for each $s\not\in E$. Set $g(s;t,x):=\psi
_{s}(t):=\sum_{p,q}d_{pq}(x)s^{p}t^{q}$. We solve (\ref{exp}) for $a_{ij}(x)$
in terms of $d_{pq}(x)$, and set $f(y,t,x)=\sum a_{ij}(x)y^{i}t^{j}$. Then, $%
f(\varphi (s,t,x),t,x)=g(s;t,x)$. Therefore, $f(y,t,x)$ diverges, and $%
\limfunc{Conv}_{\varphi }(f)=E$.
\end{proof}

\begin{corollary}
For any $f(y,t,x)\in \mathbb{C}[[y,t,x]],\ $ either $\mathfrak{c}(\limfunc{%
Conv}_{\varphi }(f))=0$ or $\limfunc{Conv}_{\varphi }(f)=\mathbb{C}$.
\end{corollary}

Finally, we list some properties of $\varphi $-convergence sets that follow
directly from the corresponding properties of polar sets.

\begin{corollary}
If $f\in \mathbb{C}[[y,t,x]]$ is divergent, then for any $F_{\sigma }$
subset $E\subseteq \limfunc{Conv}_{\varphi }(f)$, there exists a divergent
series $h(y,t,x)\in \mathbb{C}[[y,t,x]]$ such that $E=\limfunc{Conv}%
_{\varphi }(h).$
\end{corollary}

\begin{corollary}
For any finite or infinite sequence $f_{1},f_{2},...f_{n},...$ in $\mathbb{C}%
[[y,t,x]]$, then there exists an $h\in \mathbb{C}[[y,t,x]]$ such that $%
\limfunc{Conv}_{\varphi }(h)=\cup _{i=1}^{\infty }\limfunc{Conv}_{\varphi
}(f_{i}).$
\end{corollary}

\begin{corollary}
Let $f\in \mathbb{C}[[y,t,x]]$, $\lambda \in \mathbb{C}$ and $P:\mathbb{C}%
\rightarrow \mathbb{C}$ be a monic polynomial function. There are series $g$
and $h$ in $\mathbb{C}[[y,t,x]]$ such that $\limfunc{Conv}_{\varphi }(g)=$ $%
P^{-1}\left( \limfunc{Conv}_{\varphi }(f)\right) $ and $\limfunc{Conv}%
_{\varphi }(h)=\lambda \cdot \limfunc{Conv}_{\varphi }(f).$
\end{corollary}

\begin{proof}
The corollary follows since

$\mathfrak{c}\left( P^{-1}\left( \limfunc{Conv}_{\varphi }(f)\right) \right)
=\left( \mathfrak{c}(\limfunc{Conv}_{\varphi }(f)\})\right) ^{1/\deg P}$ and

$\mathfrak{c}\left( \lambda \cdot \limfunc{Conv}_{\varphi }(f)\right)
=|\lambda |\cdot \mathfrak{c}(\limfunc{Conv}_{\varphi }(f)).$
\end{proof}

Since a locally polar set is polar by Josefson's theorem (see \cite{Kl}), we
have the following.

\begin{corollary}
An $F_{\sigma }$ subset $S$ of $\mathbb{C}$ is locally a $\varphi $%
-convergence set of divergent power series if and only if it is $\varphi $%
-convergence set of a divergent power series. That is, if for every $s\in S$
there exists a neighborhood $U$ of $s$ in $\mathbb{C}$ and a divergent power
series $f_{U}$ such that $S\cap U=\limfunc{Conv}_{\varphi }(f_{U})$ then
there exists a divergent power series $f$ such that $S=\limfunc{Conv}%
_{\varphi }(f).$
\end{corollary}

\begin{example}
\emph{(cf. \cite{AM}) Let }$\varphi _{k}(t,x)$\emph{\ be a sequence of
convergent power series with }$\varphi _{k}(0,0)=0.$\emph{\ Then}%
\begin{equation}
f(y,t,x):=\sum_{n=1}^{\infty }n^{n}\prod_{j=1}^{n}\left[ y-\varphi _{j}(t,x)%
\right] \text{ \emph{and}}  \label{f}
\end{equation}%
\begin{equation}
g(y,t,x):=\sum_{i=1}^{\infty }\left[ i!t^{i}+i!^{2}(y-\varphi _{j}(t,x))^{i}%
\text{ }\right]  \label{g}
\end{equation}%
\emph{are divergent power series but for each }$k=1,2,...,$\emph{\ the
series }$f(\varphi _{k}(t,x),t,x)$\emph{\ converges while }$g(\varphi
_{k}(t,x),t,x)$\emph{\ diverges . \ The empty set and finite or countable
sets being sets of capacity zero are }$\varphi $\emph{-convergence sets of
divergent series. If }$\{s_{j}\}_{j=1}^{\infty }$\emph{\ is a sequence of
complex numbers by taking }$\varphi _{j}(t,x)=s_{j}t+\varphi (t,x)$\emph{\
for all }$j$\emph{\ in }$g(y,t,x)$\emph{\ in (\ref{g}), and by setting }$%
\varphi _{j}(t,x)=\varphi (t,x)$\emph{\ for all }$j$\emph{\ in }$f(y,t,x)$%
\emph{\ we obtain concrete examples of series with }$\limfunc{Conv}_{\varphi
}(g)=\{s_{j}\}$\emph{\ and }$\limfunc{Conv}_{\varphi }(f)=\phi .$
\end{example}

\begin{example}
\emph{The Cantor set }$C$\emph{\ in }$%
\mathbb{R}
$\emph{\ has Hausdorff dimension }$\ln 2/\ln 3,$\emph{\ hence of positive
capacity, any series }$f(y,t,x)$\emph{\ for which }$f(st+\varphi (t,x),t,x)$%
\emph{\ converges for all }$s\in C$\emph{\ is necessarily convergent.
However, there exist modifications of Cantor sets that have capacity zero,
and therefore are }$\varphi $\emph{-convergence sets of divergent power
series.}
\end{example}

\begin{remark}
\emph{The closure of a }$\varphi $\emph{-convergence set is not necessarily
a }$\varphi $\emph{-convergence set. For example, the countable set }$%
\mathbb{Q}
\ $\emph{is a }$\varphi $\emph{-convergence set of divergent }$f$\emph{\ in (%
\ref{f}) but its closure }$%
\mathbb{R}
$\emph{\ being nonpolar can not be a }$\varphi $\emph{-convergence set of a
divergent series.}
\end{remark}

\begin{remark}
\emph{The situation is quite different, as one would expect, when
restrictions of functions are considered. For example,} \emph{for any
positive integer }$k,$\emph{\ it is elementary to construct a function }$f:%
\mathbb{R}
^{n}\rightarrow 
\mathbb{R}
$ \emph{\ that is exactly }$k$\emph{-times differentiable but whose
restriction to every line in }$%
\mathbb{R}
^{n}$\emph{\ is real-analytic. The function }%
\begin{equation*}
f(x,y):=\left\{ 
\begin{array}{c}
(x^{2}+y^{2})\exp \left( -\frac{y^{2}}{x^{2}+y^{4}}\right) \text{ if }%
(x,y)\neq (0,0) \\ 
0\text{ if }(x,y)=(0,0)%
\end{array}%
\right. 
\end{equation*}%
\emph{is in }$C^{\infty }\left( 
\mathbb{R}
^{2}\right) $ \emph{\ and for all }$m\neq 0$ \emph{the single-variable
function }$f(mt,t)$\emph{\ is\ real-analytic. However, }$f$\emph{\ is not a
real-analytic function as it fails to be real-analytic along the }$y$\emph{%
-axis. Is there a }$C^{\infty }$\emph{\ function }$f:%
\mathbb{R}
^{n}\rightarrow 
\mathbb{R}
$\emph{\ which is not real-analytic but whose restriction to }every\emph{\
line is real-analytic? The answer is negative as it was shown by J. Bochnak 
\cite{bochnak} and J. Siciak \cite{siciak} that if a }$C^{\infty }$\emph{\
function }$f:%
\mathbb{R}
^{n}\rightarrow 
\mathbb{R}
$\emph{\ is real-analytic on every line segment through a point }$x_{0},$%
\emph{\ then }$f$\emph{\ is real-analytic in a neighborhood of }$x_{0}$\emph{%
. (See \cite{Ne} and \cite{Ne1} for }$C^{\infty }$-\emph{analogs of the
Bochnak-Siciak theorem.)}
\end{remark}

\end{document}